\documentclass[conference]{IEEEtran}
\IEEEoverridecommandlockouts
\usepackage{cite}
\usepackage{amsmath,amssymb,amsfonts}
\usepackage{xurl}             
\usepackage{hyperref}         
\usepackage{algorithmic}
\usepackage{graphicx}
\usepackage{textcomp}
\usepackage{xcolor}
\makeatletter 

\newcommand{\linebreakand}{%
  \end{@IEEEauthorhalign}
  \hfill\mbox{}\par
  \mbox{}\hfill\begin{@IEEEauthorhalign}
}
\makeatother 

\def\BibTeX{{\rm B\kern-.05em{\sc i\kern-.025em b}\kern-.08em
    T\kern-.1667em\lower.7ex\hbox{E}\kern-.125emX}}
\begin{document}
\title{Hardware Trends Impacting Floating-Point Computations In Scientific Applications}

\author{

\IEEEauthorblockN{Jack Dongarra}
\IEEEauthorblockA{\textit{University of Tennessee} \\
\textit{Oak Ridge National Laboratory}\\
\textit{University of Manchester}\\
Oak Ridge, TN, USA \\
dongarra@icl.utk.edu}
\and
\IEEEauthorblockN{John Gunnels}
\IEEEauthorblockA{\textit{NVIDIA Corporation} \\
Santa Clara, CA, USA \\
jgunnels@nvidia.com}
\and
\IEEEauthorblockN{Harun Bayraktar}
\IEEEauthorblockA{\textit{NVIDIA Corporation} \\
Santa Clara, CA, USA \\
hbayraktar@nvidia.com}
\and
\linebreakand
\IEEEauthorblockN{Azzam Haidar}
\IEEEauthorblockA{\textit{NVIDIA Corporation} \\
Santa Clara, CA, USA \\
ahaidarahmad@nvidia.com}
\and
\IEEEauthorblockN{Dan Ernst}
\IEEEauthorblockA{\textit{NVIDIA Corporation} \\
Santa Clara, CA, USA \\
dane@nvidia.com}
}

\maketitle
\thispagestyle{plain} 
\pagestyle{plain}    
\begin{abstract}

The evolution of floating-point computation has been shaped by algorithmic advancements, architectural innovations, and the increasing computational demands of modern technologies, such as artificial intelligence (AI) and high-performance computing (HPC). This paper examines the historical progression of floating-point computation in scientific applications and contextualizes recent trends driven by AI, particularly the adoption of reduced-precision floating-point types. The challenges posed by these trends, including the trade-offs between performance, efficiency, and precision, are discussed, as are innovations in mixed-precision computing and emulation algorithms that offer solutions to these challenges. This paper also explores architectural shifts, including the role of specialized and general-purpose hardware, and how these trends will influence future advancements in scientific computing, energy efficiency, and system design.
\end{abstract}


\section{Introduction}
Floating-point computation is foundational to modern scientific applications, enabling the representation of real numbers across a wide range of magnitudes and providing the precision necessary for calculations in fields like physics, chemistry, and engineering. Over the decades, the evolution of floating-point computation has been influenced by the increasing complexity of scientific problems, technological advancements, and the rise of new computational paradigms, such as deep neural network- (DNN-) based AI algorithms~\cite{lecun2015deep}.

\begin{figure}
\centering
\includegraphics[width=1.0\linewidth]{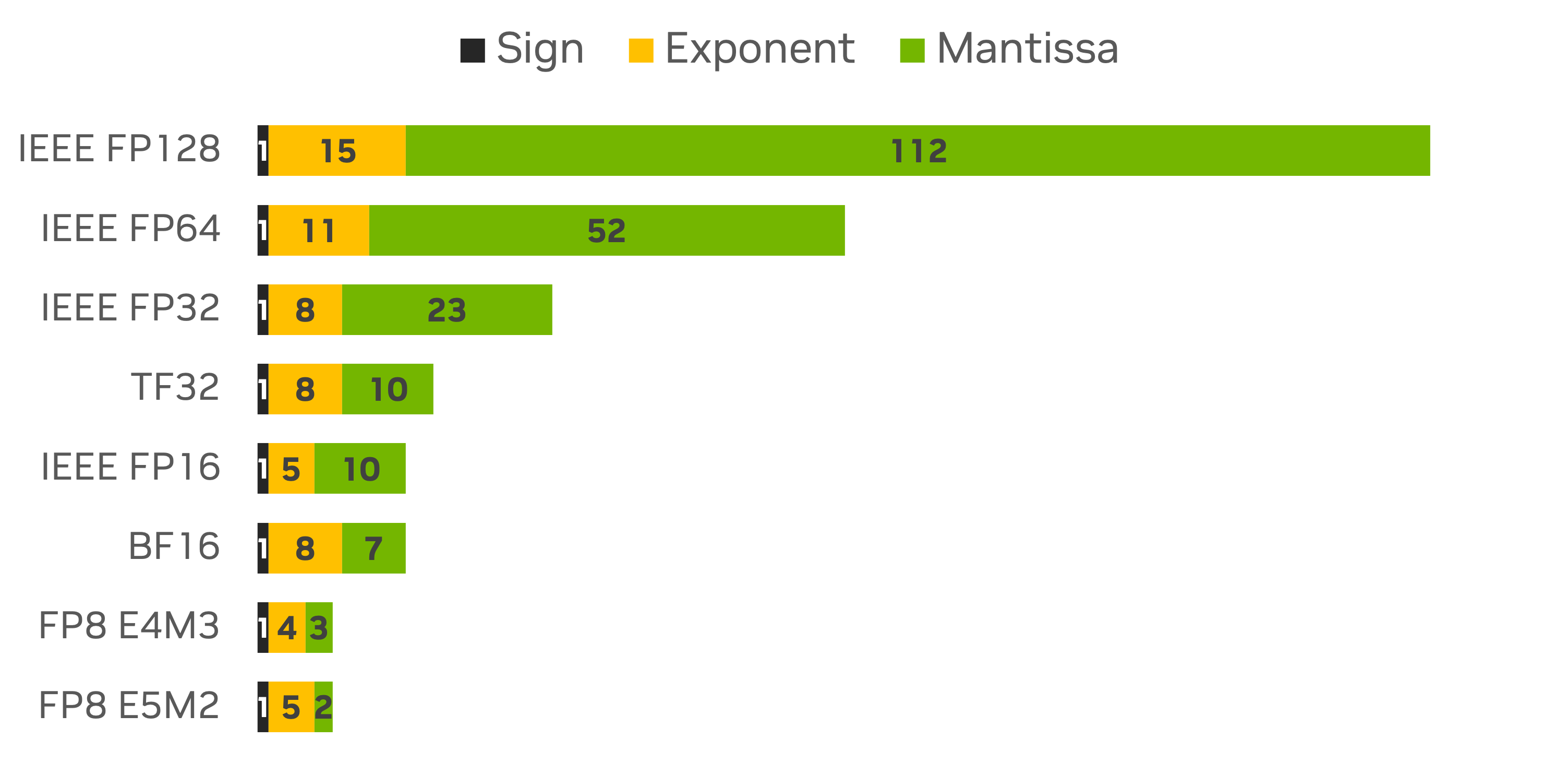}
\caption{Various floating-point (FP) representations used today in scientific computing and AI. The exponent bits determine the dynamic range of the FP number, while the mantissa bits determine the precision. Four IEEE FP types are shown: half (FP16), single (FP32), double (FP64), and quad (FP128). TensorFloat-32 (TF32), available on NVIDIA GPUs starting with the Ampere architecture, is a Tensor Core matrix multiply compute mode where input and output are FP32, but input operands are truncated. Bfloat16 (BF16), which was introduced by Google~\cite{kalamkar2019studybfloat16deeplearning}, has the same range as FP32 at the expense of mantissa bits. Two variants of FP8, with different splits of exponent and mantissa bits~\cite{micikevicius2022fp8formatsdeeplearning} are shown.}\label{fig:fp_types}
\end{figure}

This paper explores the history of floating-point computation, focusing on architectural innovations that have shaped the current landscape. It examines key developments, from early emulation to dedicated hardware, and highlights recent trends, including mixed-precision computing and reduced-precision floating-point types (Figure~\ref{fig:fp_types}). The impact of these trends on scientific computing and AI is analyzed, along with the challenges they present regarding system design, energy efficiency, and performance.

\section{Benchmarks} \label{sec:benchmarks}
Throughout this paper, we will refer to several community benchmarks that have emerged over time, each serving a critical role in evaluating and exposing performance characteristics and limitations of the underlying hardware, as well as providing a readily conveyed, widely understood record of progress. These benchmarks have become standard tools in the high-performance computing (HPC) community for assessing the efficiency and effectiveness of various computing systems. While in this paper, we focus on floating-point operation-focused benchmarks, other benchmarks exist. One such example is the Graph 500~\cite{graph500}, which measures a system's performance on graph-based problems important for large dataset analysis.

The most well-known benchmark is HPL~\cite{hpl} (High-Performance Linpack). Traditionally used to rank systems in the TOP500 list~\cite{Dongarra2011,Top500}, which focuses on a system's ability to solve dense linear equations, it measures the floating-point computing power of supercomputers, highlighting their raw computational capability. As valuable as it is, HPL emphasizes peak performance of double-precision dense matrix multiplications, which may not always correlate with real-world application performance.

In contrast, the Green500~\cite{green500} list focuses on the energy efficiency of the HPL benchmark, measuring the FLOPS per watt delivered by a system. As power consumption becomes an increasingly critical factor in supercomputing, with recent systems approaching 40 megawatts of power consumption~\cite{Top500}, Green500 plays a pivotal role in pushing the development of energy-efficient architectures and balancing performance with power usage.

HPCG~\cite{hpcg} (High-Performance Conjugate Gradient) was introduced to provide a more comprehensive measure of real-world application performance, especially for systems that perform well in terms of memory bandwidth, network latency, and irregular memory access patterns. HPCG aims to capture a broader range of system performance characteristics, some of which HPL might overlook, providing insight into a machine's ability to handle more complex, memory-bound workloads.

A more recent addition, HPL-MxP (formerly called HPL-AI)~\cite{hplmxp}, or HPL for mixed-precision, is designed to benchmark systems optimized for AI and machine learning workloads. By focusing on mixed-precision operations, HPL-MxP reflects the growing need for systems capable of efficiently handling lower-precision calculations typical in AI models, exposing the performance capabilities of modern hardware in these emerging domains without abandoning the needs of scientific computing.
\begin{figure}
    \centering
    \includegraphics[width=1.0\linewidth]{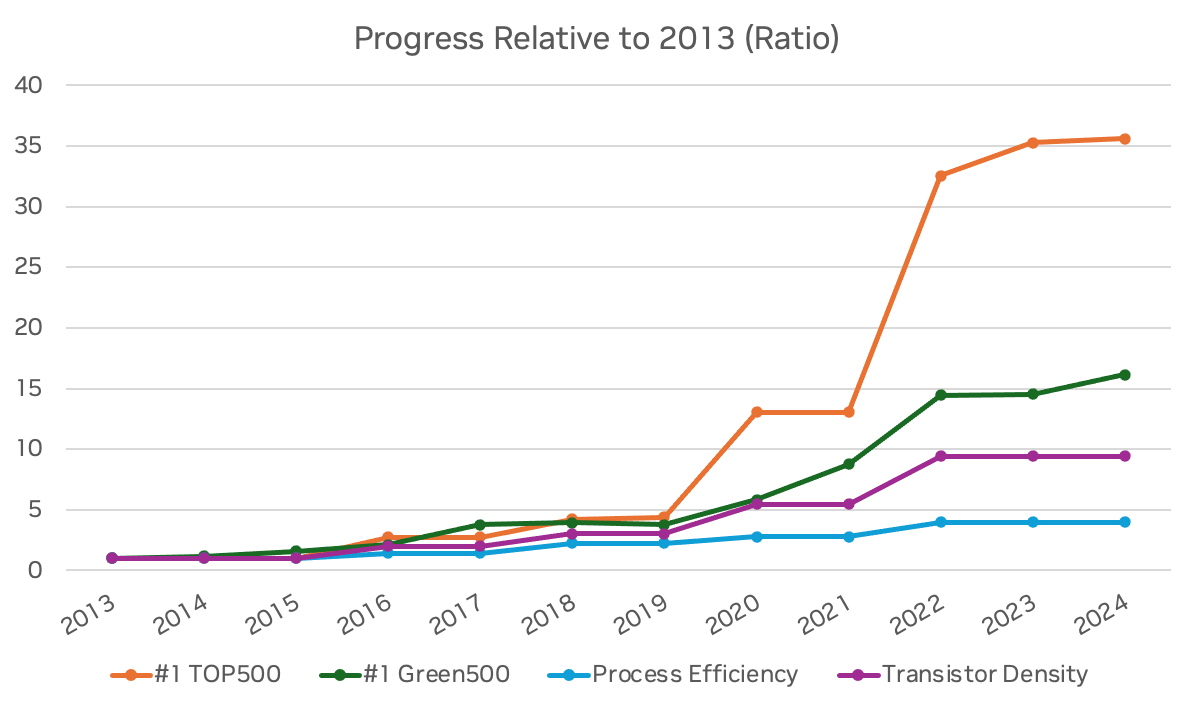}
    \caption{Historical record of the advances made in transistor density and process efficiency~\cite{sjprocessdataieee}, contrasted with the increases seen in the TOP500 and Green500 lists since 2013.  All series are normalized to 1.0 at the outset of the chart in 2013.  It is notable that both the TOP500 and Green500 entries have improved at a far greater rate than process technology and, further, that the TOP500 (performance) has increased at a greater rate than the Green500 (efficiency).}
    \label{fig:hpl-greenprocess}
\end{figure}

\begin{figure}
    \centering
    \includegraphics[width=1.0\linewidth]{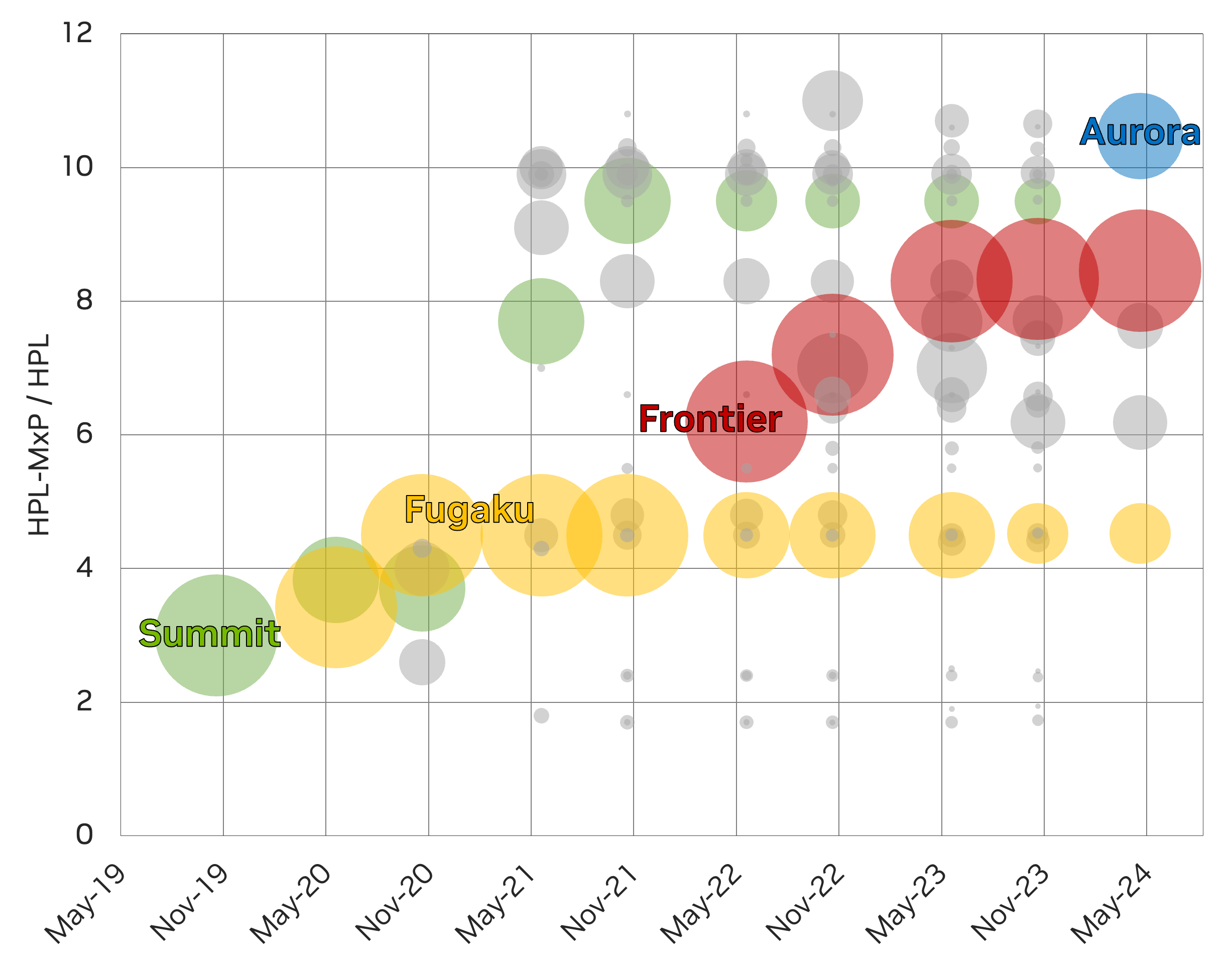}
    \caption{Ratio of HPL-MxP (formerly HPL-AI) $R_{max}$ to HPL $R_{max}$ over time, since the inception of the HPL-MxP benchmark. Some top-ranked Supercomputers are highlighted with colors and labels. The bubble size is inversely proportional to the Top500 HPL ranking of that particular supercomputer, which explains why Summit or Fugaku bubbles shrink over time. In addition to the general trend of the ratio increasing over time for top-ranked systems, implementation optimizations also increase the speed-up over time, as illustrated by the Summit and Fontier systems.}
    \label{fig:hpl-mxp}
\end{figure}

Each of these benchmarks exposes different performance aspects of a computer system, from raw computational power and energy efficiency to memory and network bandwidth and latency to flexibility concerning floating-point precision. They give a holistic view of how well a machine will likely perform across various real-world tasks. Viewed against the landscape of history, as seen in Figures~\ref{fig:hpl-greenprocess}~and~\ref{fig:hpl-mxp}, these benchmarks can provide not only a comparator to larger hardware trends (such as processor efficiency or transistor density), but an indication of the factors that will constrain progress in future systems. Taken together, they are essential tools for guiding hardware design and system optimizations in the quest for faster, more efficient, and versatile computing platforms.\\

\section{The Evolution of Floating-Point Computation}

\subsection{Early Emulation}
In the early days of computing, into the mid-1950s, floating-point operations were typically performed via software emulation, where general-purpose processors simulated floating-point arithmetic utilizing fixed-point representations and operations~\cite{Backus:1954:ISS}. This approach, while functional, was slow and resource-intensive. Emulation required many CPU cycles for each floating-point operation, making these computations much slower than integer arithmetic.

Despite its limitations, emulation allowed early computers to perform scientific calculations and helped lay the foundation for future advancements in floating-point hardware throughout the 1960s and 1970s.

\subsection{The Co-Processor Era}
The introduction of dedicated floating-point co-processors in the 1980s marked a significant leap forward in floating-point computation. These co-processors, such as Intel’s 8087~\cite{Palmer80}, were separate hardware components designed to handle floating-point operations independently of the CPU. This separation significantly improved performance and allowed computers to tackle more complex scientific problems.

While the 8087 might be the best-known example of this technology, co-processors were in widespread use at the time.  For example, the Motorola 68020, relied on external floating-point units (FPUs), such as the 68881~\cite{Huntsman:1983:MFP}, similar to the way in which Intel's pre-x486 and AMD's pre-K5 architectures, used external units like the 8087. External FPUs were also seen in early SPARC systems and MIPS processors, such as the SPARCstation 1 with the Weitek 3170 FPU~\cite{Birman:1990:DWS} and the MIPS R4000's CP1~\cite{10.5555/154056}.

However, using co-processors introduced additional system complexity, requiring specialized hardware and coordination between the CPU and the FPU. Despite this, co-processors became a staple in high-performance systems, enabling faster scientific computing and simulations.

\subsection{Integration of Floating-Point Units}

The release of Intel’s x486~\cite{522958} processor in 1989 marked a turning point in floating-point computation. The x486 integrated the FPU directly into the CPU, eliminating the need for separate hardware. This integration simplified system design, improved performance, and made floating-point operations a standard feature of general-purpose computing.

As computing requirements expanded, especially in fields such as scientific computing and multimedia, FPUs became standard across many architectures.  With the x486, floating-point computation became more accessible and widely used in applications such as computer graphics, simulations, and scientific calculations. This development set the stage for the widespread adoption of floating-point hardware in both consumer and professional computing environments.

Beyond Intel’s x86 line, many CPU architectures incorporated FPUs, either as external coprocessors or integrated directly into the CPU, to enhance floating-point performance. Notable examples include the Motorola 68040~\cite{46770}, which was the first in the 68000 series to integrate the FPU, and the PowerPC 601~\cite{5389550}, which became popular in Apple's early Macintosh systems for its integrated floating-point capabilities.

\subsection{The GPU Revolution}

GPUs (Graphics Processing Units) were initially developed in the 1980s and 1990s to meet the growing demand for 2D and 3D graphics rendering. Early GPUs, produced by companies such as SGI (Silicon Graphics), 3dfx, NVIDIA, and ATI (later acquired by AMD), were primarily focused on enhancing the real-time rendering of images, textures, and geometry, especially in gaming and graphical user interfaces (GUIs). A significant milestone came in 1999 with the release of NVIDIA’s GeForce 256~\cite{9623445}, marketed as the first ``GPU" capable of processing graphics independently from the CPU. This allowed the CPU to focus on other tasks while the GPU specialized in real-time rendering. GPUs handled critical graphics tasks such as vertex transformations, lighting calculations, and texture mapping, all essential for 3D graphics in gaming and multimedia.

The early 2000s saw another major shift with the introduction of programmable GPUs. NVIDIA’s GeForce3, launched in 2001, was a key milestone in this revolution, offering programmable shaders that allowed developers to directly program the GPU for custom operations, including floating-point computations~\cite{buck-brook}.

The mid-2000s marked a turning point for GPUs, as their parallel processing capabilities were broadly recognized as useful in scientific and engineering computing. Unlike CPUs, which are optimized for serial processing and excel at handling a few threads quickly, GPUs have hundreds or even thousands of cores optimized for parallel workloads, making them ideal for tasks that can be divided into many smaller operations that thousands of threads can execute. This made GPUs especially useful for HPC, where they could accelerate simulations and large-scale mathematical computations, and, later, for machine learning and AI, where they became indispensable for training deep learning models. 
\begin{table}[ht]
\centering
\caption{GPU Generations: Compute Throughput Vs Memory Bandwidth}
\begin{tabular}{|l|c|c|c|c|}
\hline
\textbf{Figure of Merit} & \textbf{Volta} & \textbf{Ampere} & \textbf{Hopper} & \textbf{Blackwell}\\
\textbf{} & \textbf{(V100)} & \textbf{(A100)} & \textbf{(H200)} & \textbf{(B200)}\\
\hline
\textbf{FP64 FMA (TFLOP/s)}& 7.8  & 9.75& 33.5& 40\\
 \textbf{FP64 Tensor (TFLOP/s)}& N/A & 19.5 & 67 &40\\
\textbf{FP16 FMA (TFLOP/s)}& 31.4& 78& 134& 80\\
 \textbf{FP16 Tensor (TFLOP/s)}& 125 & 312 & 989 &2250\\
 \hline
\textbf{Memory BW (TB/s)}& 0.9& 2.0& 4.8& 8.0\\
\hline
\textbf{FP64 FMA (B/FLOP)}& 0.124& 0.225& 0.158& 0.220
\\
 \textbf{FP64 Tensor (B/FLOP)}& N/A & 0.112& 0.079&0.220
\\
\textbf{FP16 FMA (B/FLOP)}& 0.031& 0.028& 0.039& 0.110
\\
 \textbf{FP16 Tensor (B/FLOP)}& 0.008& 0.007& 0.005&0.004\\
 \hline
\end{tabular}
\label{table:hw_specs}
\end{table}

NVIDIA’s CUDA~\cite{NickollsBGS08} (Compute Unified Device Architecture) platform, first introduced in 2006, enabled a broad developer base to harness GPU power for non-graphics workloads. As a result, GPUs, initially designed for rendering graphics, soon became indispensable for general-purpose computing tasks, particularly in scientific computing~\cite{4541126, liu2008accelerating, GPUAccelerated}. Breakthroughs in AI algorithms came soon after~\cite{raina2009large, NIPS2012_c399862d} and rapidly led to leveraging GPUs outside of scientific computing. The GPU's ability to perform parallel floating-point operations at high speeds, complemented by equally impressive bandwidth capabilities (see Table~\ref{table:hw_specs}), firmly established them as ideal for diverse tasks requiring massive computational power, from training deep learning models to running large-scale simulations in physics and biology.

Today, GPUs are widely used in fields as varied as engineering and scientific simulations, medical research, cryptocurrency mining, big data processing, finance, and countless AI applications, cementing their role as a critical tool for a wide range of performance-critical computational tasks beyond graphics rendering.\\

\section{Recent Trends in Floating-Point Computation}

\subsection{The Rise of AI and Reduced Precision}
Artificial intelligence, particularly deep learning, has profoundly impacted floating-point computation. AI workloads are dominated by matrix multiplications and tensor operations which can tolerate lower precision without a significant loss in accuracy. As a result, reduced-precision floating-point types, such as FP16 (half-precision), BF16 (brain floating-point, also referred to as bfloat16), and, more recently, FP8~\cite{micikevicius2022fp8formatsdeeplearning, ocp_fp8, rouhani2023microscalingdataformatsdeep} have become increasingly popular in AI applications (see Figure~\ref{fig:fp_types}).

Reduced precision offers several advantages, including faster computations and lower energy consumption. In AI, especially in training and inference tasks, models can maintain high accuracy using lower precision, which leads to improved throughput and efficiency. This shift has driven the development of specialized hardware, such as NVIDIA’s Tensor Cores~\cite{nvidia2017v100}, initially optimized for matrix operations at a reduced precision.

It should be noted that the move toward reduced-precision capabilities poses challenges for applications that require high accuracy, such as scientific simulations and financial modeling. While AI applications can often tolerate lower precision, many scientific fields depend on high-precision calculations to ensure the validity of their results.

\subsection{Mixed-Precision Computing}
Mixed-precision computing has emerged as a promising solution to address the limitations of reduced precision. Mixed-precision algorithms~\cite{baboulin2009accelerating}, distinct from the mixed-precision operations they sometimes leverage, dynamically adjust the precision of floating-point calculations based on the accuracy required for each specific task. This approach allows systems to use lower precision for less critical calculations while reserving higher precision for tasks that require greater accuracy (see~\cite{mixed}), requiring careful algorithmic design to ensure that performance gains from using lower precision do not come at the expense of accuracy. 

In the discussion of mixed-precision computation that follows, an approach is explored that involves using different levels of precision at various stages of the computation process for solving a dense system of equations. Specifically, one phase of the computation, where the highest level of accuracy is not critical, may be performed using a lower (or even fixed) precision, such as half-precision (FP16), to optimize performance and reduce resource consumption~\cite{haidar2018design}.

As the computation progresses, the algorithm switches to a higher precision, such as single-precision (FP32) or double-precision (FP64), to complete parts of the process where increased accuracy is essential. This transition is necessary in phases where errors accumulated during earlier steps must be corrected or refined, ensuring that the final result meets the desired accuracy thresholds. The dynamic adjustment between lower and higher precision enables a balance between computational speed, energy efficiency, and numerical precision. By strategically employing mixed precision, we can achieve significant performance gains without compromising the overall accuracy of the computation.

Mixed-precision computing has proven especially useful in AI~\cite{pmlr-v37-gupta15} and HPC~\cite{CLARK20101517}. For example, in deep learning, models can be trained using a combination of FP8, BF16, FP16, and FP32 calculations, reducing the computational load without sacrificing accuracy.  The same domain affords opportunities for finer-grained application of this technique. For example, forward propagation through a neural network typically requires higher precision for weights and activations. In contrast, gradients in the backward propagation (used for updating weights) require a higher dynamic range. This has led to the introduction of two different FP8 variants, E4M3 and E5M2~\cite{micikevicius2022fp8formatsdeeplearning}. Mixed-precision algorithms must manage these precision transitions efficiently to maximize performance without forfeiting the quality of results.
\begin{table}[hbt]
    \centering
    \caption{Mixed-precision iterative refinement solver (from cuSOLVER) performance and efficiency for solution of a 32k double-precision complex system of equations}
    \begin{tabular}{|l|l|c|c|c|} 
    \hline
         &  &  \textbf{Volta} &  \textbf{Ampere} & \textbf{Hopper} \\
         &  &  \textbf{(V100)} &  \textbf{(A100)} & \textbf{(H200)} \\
         \hline
         \hline
         Performance& FP64&    6.6&  16.9& 42.6\\
         TFLOP/s& FP16+FP64 MxP&    34.6&  74.6& 124.2\\
         \hline
         Efficiency& FP64&    28&  45& 78\\
         GFLOP/s/Watt& FP16+FP64 MxP&    173&  262& 529\\
         \hline
    \end{tabular}
    \label{tab:mxp_perf}
\end{table}

\begin{figure}[hbt!]
\centering
\includegraphics[width=1.0\linewidth]{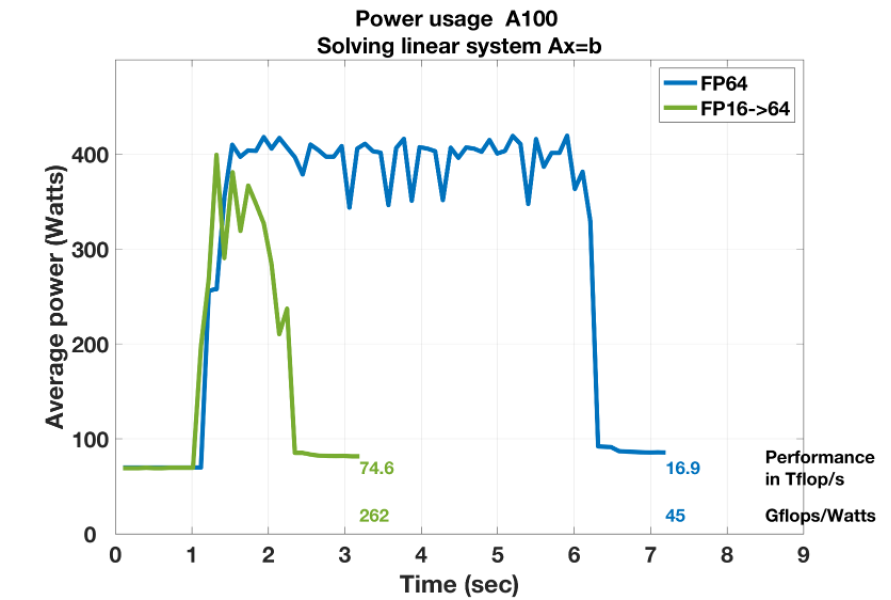}
    \caption{Representative power consumption curves measured on an NVIDIA Ampere A100 GPU during the execution of two different equations solvers is shown. The blue line shows the FP64 LU solver (corresponds to ZGETRF \& ZGETRS in LAPACK) while the green line shows the Tensor Core accelerated mixed-precision iterative refinement solver available in the cuSOLVER library, which relies on cuBLAS for Level 3 BLAS operations, both for a matrix size of 32000. Comparison with other GPU architectures can be found in Table~\ref{tab:mxp_perf}.}
    \label{fig:tcairs}
\end{figure}

This technique has also been applied in scientific computing, in cases where certain phases of simulations can tolerate lower precision, allowing for faster computations (see Table~\ref{tab:mxp_perf} and Figure~\ref{fig:tcairs}).  Using error-correction techniques, such as stochastic rounding~\cite{croci2022stochastic} and iterative refinement to prevent the propagation of errors, allows mixed-precision algorithms to maintain high accuracy even when using reduced precision for certain operations, making mixed-precision computing suitable for improving the performance of various scientific and AI applications.

By optimizing precision based on task requirements, mixed-precision computing also improves energy efficiency, which is a critical constraint for modern floating-point computation. Figure~\ref{fig:tcairs} shows an example of a mixed-precision solver that is 4.4 times faster and 5.8 times more power-efficient than one using full precision for the entire algorithm.

\subsection{Emulation in Modern Systems}
Emulation, once a necessity in the absence of dedicated floating-point hardware, has come back in modern systems.  Emulation techniques have evolved in terms of their underlying methods, implementations, and the flexibility with which they leverage (plentiful) hardware resources~\cite{ooyo22, uchino2024performanceenhancementozakischeme, henry2019leveragingbfloat16artificialintelligence}. This adaptability is particularly important in domains such as scientific AI and HPC, where computational needs often exceed the capabilities of the more established, higher-precision hardware resources.
\begin{table}[hbt!]
    \centering
    \caption{Preliminary HPL Performance and Efficiency Measurements on Blackwell B200 GPU comparing emulation with s=7 (INT8 data elements)~\cite{uchino2024performanceenhancementozakischeme} to native FP64 [data subject to change]}
    \begin{tabular}{|l|l|c|c|c|} 
    \hline
         &  &  \textbf{Native}&  \textbf{Emulation}& \\
         &  &  \textbf{FP64}&  \textbf{(s=7)}& \textbf{Ratio}\\
         \hline
         \hline
         At Maximum& TFLOP/s&    34.5&  68.4& 2.0\\
         Performance& GFLOP/s/Watt&    41.7&  71.3& 1.7\\
         \hline
         At Maximum& TFLOP/s&    23.1&  53.4& 2.3\\
         Efficiency& GFLOP/s/Watt&    51.4&  82.1& 1.6\\
         \hline
    \end{tabular}
    \label{tab:emulation_perf}
\end{table}

For example, high-precision floating-point operations can be emulated on hardware designed for lower, or even fixed-, precision, allowing systems to balance performance and accuracy. Emulation also provides a novel means to support evolving computational demands without requiring a complete hardware overhaul and something of an added degree of freedom in designing new systems, as it allows some components to serve a dual purpose. This is demonstrated in Table~\ref{tab:emulation_perf} where the performance of HPL employing native FP64 and emulation using the techniques described in~\cite{uchino2024performanceenhancementozakischeme} are compared. When configured for maximum performance, emulation yields a two-fold improvement in performance and 70\% improvement in power efficiency. In contrast, when maximum efficiency is the target, emulation yields a 2.3x speed-up and a 60\% improvement in power efficiency.\\

\section{Challenges in Floating-Point Computation}

\subsection{Balancing Precision and Efficiency}
One of the key challenges in floating-point computation is finding the right balance between precision and efficiency. High-precision formats like FP64 (double-precision) are necessary for certain scientific applications, but they come at a cost regarding speed and energy consumption, as seen in Figure~\ref{fig:hpl-mxp}. In contrast, reduced-precision formats like FP16 are much faster and more efficient but may not provide the accuracy needed for all tasks (See Figure \ref{fig:fp_types}).

This trade-off is particularly pronounced when considering fields like AI, where reduced precision is sufficient for many tasks, and scientific computing, where the accuracy of results cannot be sacrificed. Today, the same systems often run applications of both stripes and the integration of AI methods and scientific computing is a burgeoning field of study~\cite{lavin2022simulationintelligencenewgeneration}.  As a result, system designers must carefully consider the needs of each application when selecting the appropriate precision format.

\subsection{Specialized vs. General-Purpose Hardware}
Another challenge in floating-point computation is the tension between specialized hardware (such as custom FPUs with the ability to execute compound instructions) and general-purpose processors (CPUs). On the one hand, specialized hardware is optimized for specific tasks, such as matrix multiplications in AI, but may lack the flexibility needed for broader applications. General-purpose processors, on the other hand, can handle a wide range of tasks but are not as efficient for specialized computations.\\

\section{Architectural Innovations and Their Impact}

\subsection{Heterogeneous Computing}
Heterogeneous computing has become a cornerstone of modern floating-point computation, combining CPUs, GPUs, and other accelerators to optimize performance. By accelerating floating-point operations via different types of processors, systems can achieve higher performance and energy efficiency.

Unsurprisingly, the tolerance that programmers once had for (relatively) distant accelerator engines, in terms of both programmability and access latency, has decreased over the years and concerted efforts have been undertaken to enable practitioners to ``eat their cake and have it too.''  Almost a decade ago, the Sierra and Summit supercomputing systems, at Lawrence Livermore (LLNL) and Oak Ridge (ORNL) National Laboratories, respectively, leveraged the first generation of NVLINK to tightly couple IBM's POWER9 processor to NVIDIA's Volta (V100) GPU~\cite{10.1147/JRD.2018.2846978}.  Today, NVIDIA’s Grace Hopper Superchip architecture~\cite{gh1} integrates tightly coupled CPU and GPU components, enabling seamless transitions between general-purpose and specialized processing. This integration allows systems to handle highly serial sparse calculations, high-throughput scientific calculations, and reduced-precision AI workloads efficiently. AMD's Instinct MI300A and Apple's M3~\cite{tandon2024portinghpcapplicationsamd,keycdna,appleAppleUnveils} are additional realizations of the goal to tightly couple distinct compute resources in a way that allows them to be viewed less as separate entities and more as a potent, unified resource.

Entwined with the drive for tighter integration, the trend of designing custom silicon tailored to specific software needs has grown significantly. Companies like Google and Apple have invested in designing chips, such as Google’s TPUs for the acceleration of AI workloads and Apple’s A-series chips for sophisticated mobile devices. These chips are designed in tandem with the software they will run, allowing for optimizations that general-purpose hardware cannot achieve. This hardware-software co-design allows for significant improvements in performance and efficiency, as the hardware is fine-tuned to the software’s needs.

Addressing the needs of the developer community, heterogeneous computing environments~\cite{10.1145/3636480.3637097,10.1145/2788396,10.1145/1542476.1542525, 6468503} offer several advantages, including the ability to optimize each task for the most suitable processor. This approach not only improves performance, but also reduces energy consumption for a broad range of applications by offloading computationally intensive tasks to specialized hardware.

\subsection{Energy Efficiency in Floating-Point Hardware}
Energy efficiency has become a critical consideration in the design of floating-point hardware across a wide variety of computing environments, from large-scale data centers and supercomputers to battery-reliant mobile devices such as laptops and phones. As computational workloads continue to grow, the energy required to power these platforms has become a significant constraint, with some public supercomputers approaching 40 megawatts. The centrality of energy consumption is reflected in the growing attention given to the Green500 list, described in Section~\ref{sec:benchmarks}.

Introduced by NVIDIA with the Volta~\cite{jia2018dissectingnvidiavoltagpu} architecture in 2017, Tensor Cores are specialized hardware designed to accelerate matrix multiply-accumulate (MMA) operations, which are critical for deep learning tasks such as matrix multiplication and convolution in neural networks. These cores are particularly efficient due to the use of complex instruction types optimized for high-throughput operations (Table~\ref{table:hw_specs}). Tensor Cores are available across many different precisions and are highly utilized in AI workloads that make extensive use of dense matrix operations. This complex operation approach significantly speeds up computations and reduces energy requirements without compromising requisite accuracy (see Figure~\ref{fig:tcairs} and Table~\ref{tab:mxp_perf}).

Tensor Cores~\cite{nvidiaTensorCoresSite} excel at accelerating matrix multiplications (\(A[m \times k] \cdot B[k \times n]\)), a fundamental operation in neural networks. They process these operations in blocks (e.g., (m, n, k) = (16, 16, 8))~\cite{nvid24b}, which boosts throughput for critical tasks like forward and backpropagation. This allows neural network operations that typically require hundreds of regular GPU instructions to be executed with relatively few tensor operations executed in a fraction of the time. Furthermore, Tensor Cores enable massive parallelism by executing multiple floating-point operations simultaneously and in a systolic manner, resulting in high throughput for deep learning models such as convolutional neural networks (CNNs)~\cite{NIPS2012_c399862d} or transformers~\cite{10.5555/3295222.3295349}. Additionally, their power efficiency is notable, as the ability of a single Tensor Core to perform a chain of operations that would require multiple traditional GPU instructions reduces energy consumption—a key consideration for both chip-level performance limits and the capabilities of large-scale data centers. Thus, the development of low-power FPUs and specialized processors has enabled systems to meet the growing demand for computational power without exceeding energy constraints.\\

\section{The Role of Floating-Point in AI and HPC}
\subsection{AI-Driven Workloads}
Artificial intelligence, particularly deep learning, has revolutionized floating-point computation. AI workloads are characterized by large-scale matrix multiplications requiring massive computational power. Reduced-precision floating-point types like FP16 have become the standard for these workloads, offering the best balance between performance and accuracy~\cite{micikevicius2018mixed}.

Tensor operations, which are the foundation of most AI models, benefit from the parallel processing capabilities of GPUs and Tensor Cores. These specialized processors are optimized for matrix operations, allowing AI models to be trained and deployed more quickly and efficiently.

As AI grows in importance, the demand for floating-point hardware that can handle AI-specific workloads will increase. This trend will drive further innovations in reduced-precision computing and specialized hardware for AI.

\subsection{Scientific Computing and High Precision}
In contrast to AI, scientific computing often requires high-precision floating-point calculations. Most simulations in fields like molecular dynamics, computational mechanics, and fluid dynamics depend on FP64 precision to ensure the accuracy of their results. These calculations are typically run on HPC systems designed to handle large-scale simulations that require high precision and significant computational power.

While reduced-precision techniques are becoming more common in scientific computing~\cite{higham2019simulating, ImprovingWeatherForecastSkillthroughReducedPrecisionDataAssimilation}, they are often used in conjunction with high-precision calculations. 
This approach improves the speed and efficiency of scientific simulations, enabling researchers to run more complex models in less time without sacrificing accuracy.\\

\section{Hardware's Impact on Software}
The interaction between hardware and software is deeply interdependent, with hardware setting the constraints within which software must operate. As hardware evolves, it directly influences how software is designed, written, and optimized. Co-evolution drives both fields forward, as hardware improvements open new opportunities for software innovation while software demands push hardware development to new heights.

\subsection{Performance and Capabilities}
The performance of hardware, particularly in terms of floating-point computation, dictates the upper limits of what software can achieve. For example, the availability of GPUs and Tensor Cores with specialized floating-point capabilities has enabled popular AI frameworks like PyTorch~\cite{paszke2019pytorchimperativestylehighperformance}, TensorFlow~\cite{10.5555/3026877.3026899}, and JAX~\cite{jax2018github} to handle large-scale matrix operations more efficiently. As a result, software developers can design more complex models and algorithms that rely on the enhanced performance provided by these hardware features without undue consideration of architectural details.

\subsection{Instruction Sets and Architectures}
On CPUs, the hardware architecture, such as x86, ARM, or RISC-V, determines the instruction sets (ISAs) available to software. Software must be compatible with the hardware's architecture, which affects how low-level operations are performed and how efficiently the software can execute. As an example, the introduction of SIMD (Single Instruction, Multiple Data) and specialized floating-point instructions has allowed scientific and AI software to accelerate matrix operations and other floating-point-intensive tasks.

Similarly, on GPUs, the Parallel Thread Execution (PTX)~\cite{nvid24b} ISA is the interface to enable it as a computing device. PTX is a part of CUDA, the parallel computing platform developed by NVIDIA that enables GPUs to perform general-purpose computing tasks beyond graphics rendering. It allows developers to leverage the massive parallel processing power of GPUs by offloading compute-intensive tasks from the CPU. CUDA provides a unified architecture that manages thousands of threads in a SIMT model (Single Instruction, Multiple Threads), supports various precision modes (such as FP16 and FP64), and optimizes memory access, making it highly efficient for parallel workloads like scientific computing, AI, and machine learning.

As NVIDIA’s GPUs have evolved, CUDA has remained central to utilizing them to maximum advantage, particularly with the introduction of Tensor Core instructions for mixed-precision operations, beginning with the Volta architecture. These innovations have led to dramatic increases in compute throughput, especially in AI/ML applications (Table~\ref{table:hw_specs}). Memory bandwidth has also significantly improved, offering higher Bytes/FLOP for non-Tensor Core operations. CUDA’s ability to manage parallel execution, memory hierarchy, and precision makes it essential for extracting maximum performance from modern GPUs.

\subsection{Specialized Hardware and Software Paradigms}
Developing specialized hardware, such as GPU Tensor Cores and Tensor Processing Units (TPUs)~\cite{10.1145/3140659.3080246}, has created new software paradigms. For instance, deep learning frameworks are optimized to leverage the parallelism of GPUs, which has drastically improved the training times for neural networks. Without these hardware advances, modern AI software could not scale to the levels required for training models like GPT-4~\cite{openai2024gpt4technicalreport} or other large neural networks.  Conversely, this is an example of the inseparability of hardware and software.  Because software has been able to leverage hardware capabilities to great advantage, technology has been pushed to deliver ever greater resources to supply the needs of applications.\\
 
\section{Software's Impact on Hardware}

Software has increasingly become a driving force in hardware design, as complex and demanding applications push the limits of existing hardware capabilities. As software grows more intricate, hardware must evolve to meet the demands for greater performance, efficiency, and flexibility.

\subsection{Energy Efficiency and Power Constraints}

As software applications become more resource-intensive, hardware must prioritize energy efficiency to maintain performance without consuming unsustainable amounts of power. Big data, real-time analytics, and AI-driven applications have all contributed to a demand for hardware that can deliver high performance per watt. In response, hardware manufacturers have developed energy-efficient architectures, such as Arm processors for data centers, GPUs with dynamic voltage and frequency scaling (DVFS)~\cite{10.5555/1924920.1924921}, as well as applications of this technology to commodity desktop and server processors, for example Intel's SpeedStep~\cite{genossar2003intel} and AMD's Cool’n’Quiet~\cite{10.1145/2208828.2208840}.\\

\section{Non-Standard Data Types}

While floating-point representations, from FP8 to FP64, dominate AI and scientific applications, several non-standard data types, and their corresponding computational mechanisms, have emerged in recent years. These data types offer marked potential advantages in terms of precision, energy efficiency, and computational speed, particularly in specialized applications.  Exotic formats and technologies such as posits~\cite{posit}, Spiking Neural Networks (SNNs)~\cite{spinnaker}, and analog computing~\cite{maley2011analog} offer tantalizing potential advantages along several axes of interest, but they face significant hurdles in terms of hardware and software support, scalability, and noise management. As computing demands evolve, these data types could find greater adoption in specialized fields where their advantages are most beneficial and their shortcomings are acceptable or less keenly felt, for example, in environments requiring extreme levels of energy efficiency or domains with specific accuracy requirements.\\

\section{Algorithmic Complexity and Memory Bandwidth}

\subsection{Impact of Algorithmic Complexity on Floating-Point Performance}
The complexity of algorithms significantly influences the efficiency of floating-point operations.  For example, dense linear algebra algorithms, such as matrix multiplications, exhibit high floating-point operation intensity (\textit{a.k.a.} arithmetic intensity) and are well-suited to GPUs with high floating-point throughput. In contrast, sparse matrix operations often require irregular memory accesses, leading to memory bandwidth bottlenecks that reduce floating-point efficiency. This impact is so broadly and acutely felt, especially in contexts that require real-time data processing, that multiple high-profile benchmarks, most notably HPCG (see section~\ref{sec:benchmarks}), provide a figure of merit for systems, largely based on this characteristic.  As a result, optimizing for memory access patterns becomes a critical aspect of hardware-software co-design in HPC and AI. Table~\ref{table:hw_specs} shows pertinent GPU specifications for FP16 and FP64 computations and the corresponding Bytes/FLOP values which are important for application performance.

\begin{figure}
\centering
\includegraphics[width=1.0\linewidth]{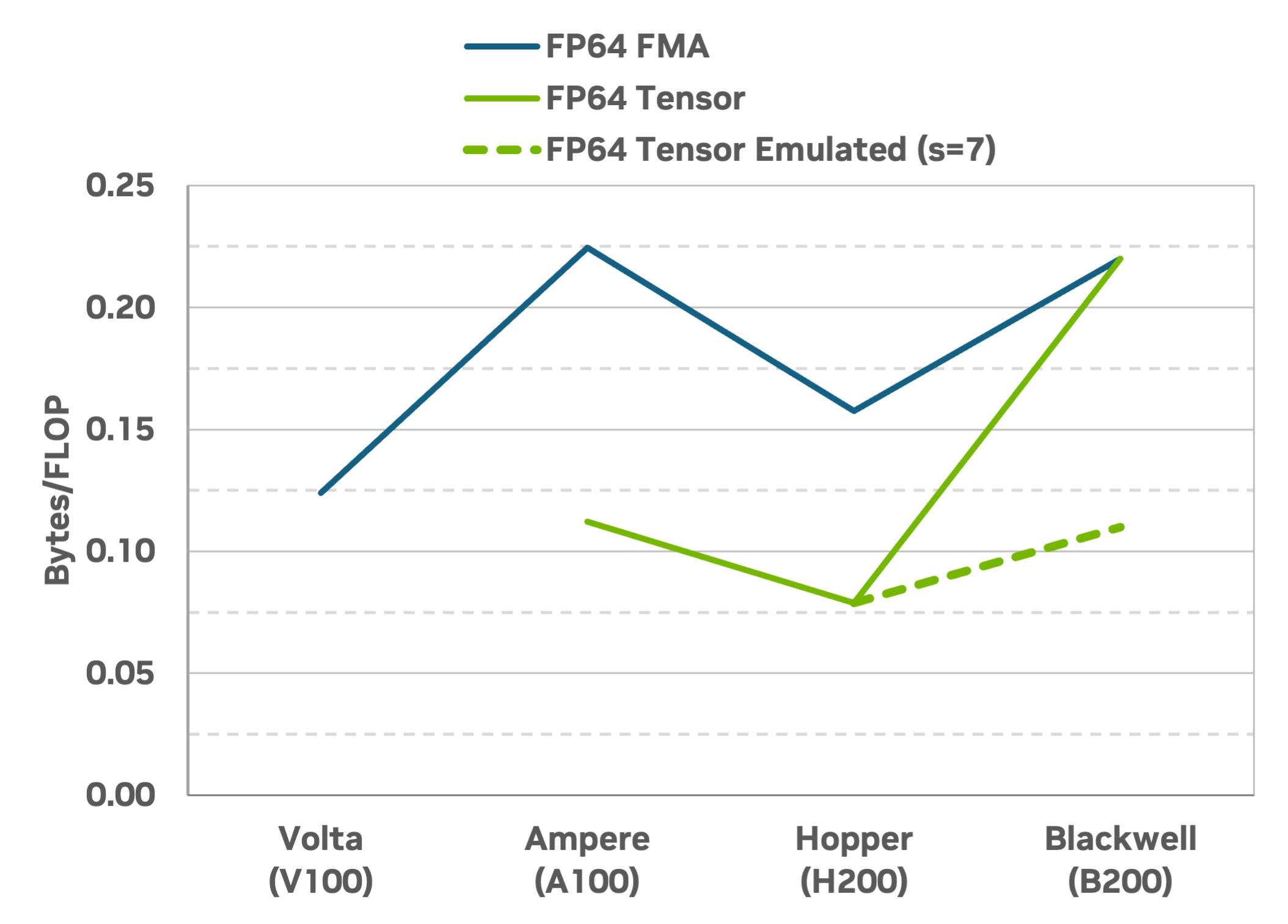}
    \caption{Comparison of Bytes/FLOP across four generations of GPUs for FMA and Tensor Core throughput from Table~\ref{table:hw_specs}. The dashed line shows Tensor Core accelerated DGEMM performance using integer based emulation with 7 slices (INT8 data storage elements)~\cite{uchino2024performanceenhancementozakischeme}. }
    \label{fig:Bytes_per_Flop}
\end{figure}

The introduction of high-bandwidth memory (HBM)~\cite{7939084} and on-chip memory hierarchies has helped mitigate some of these issues by providing faster data access and reducing the latency associated with memory operations. Figure~\ref{fig:Bytes_per_Flop} plots the data from Table~\ref{table:hw_specs}, illustrating improvements over time in the data access per floating-point operation over generations, which is important for a very broad class of applications that cannot leverage high-throughout matrix multiply operations.

Even with these hardware innovations, one of the remaining key challenges in HPC and AI applications is performance optimization across high and low arithmetic intensity components of complex algorithms. Overcoming these challenges often involves a technique called kernel fusion, where large amounts of repeated data load and store operations are avoided by combining multiple kernels, sometimes at the expense of performing some extra floating-point operations. A good example of this, from Transformers for LLM applications, is Flash Attention and its variants~\cite{dao2022flashattentionfastmemoryefficientexact,dao2023flashattention2fasterattentionbetter}.\\

\section{Future Directions in Floating-Point Computation}

\subsection{Advances in Mixed-Precision Techniques}
As AI and scientific computing evolve, mixed-precision computing will become increasingly important. Future systems will likely incorporate more sophisticated algorithms that dynamically adjust precision levels to optimize performance and energy efficiency. These systems will be able to switch between FP8, BF16, FP16, FP32, and FP64 as needed, ensuring that each task is handled with the appropriate level of precision.

This flexibility will be particularly valuable in environments where both AI and scientific computing tasks are performed, as it will allow systems to optimize for both speed and accuracy without compromise.

\subsection{Emulation for Flexibility}
Emulation will continue to play a key role in floating-point computation, particularly as new applications require higher precision or more specialized calculations. Emulation provides a flexible solution for extending the capabilities of existing hardware, allowing systems to perform floating-point operations that are not natively supported by the hardware, and offering power efficiency gains (see Table~\ref{tab:emulation_perf}).

As computational demands grow, emulation will become an increasingly valuable tool for maintaining flexibility and extending the lifespan of hardware systems.

\subsection{Energy-Efficient Designs}
The need for energy-efficient floating-point hardware will only increase as computational workloads grow. Future innovations in floating-point design will focus on reducing power consumption while maintaining high performance. This will involve the development of more energy-efficient FPUs and new architectures that minimize the energy cost of common floating-point computational patterns.
These innovations will be significant in data centers and HPC environments, where energy consumption is a major constraint on system performance.\\
 
\section{Conclusion}
The evolution of floating-point computation has been driven by advancements in both hardware and software, shaped by the demands of scientific research, artificial intelligence, and high-performance computing. As the field continues to evolve, innovations in mixed-precision computing, emulation, energy-efficient design, and non-standard data types will play a critical role in meeting the growing demands of modern applications.
By balancing the competing needs for precision, performance, and energy efficiency, floating-point hardware will remain a key component of scientific and AI-driven computation, enabling future breakthroughs in research and technology.

\bibliographystyle{ieeetr}

\vspace{12pt}
\end{document}